\documentclass[10pt]{amsart}

\usepackage{amsmath}
\usepackage{amssymb}
\usepackage[all]{xy}

\makeatletter
\@addtoreset{equation}{section}
\makeatother

\newtheorem{theorem}{Theorem}[section]
\newtheorem{proposition}[theorem]{Proposition}
\newtheorem{lemma}[theorem]{Lemma}
\newtheorem{corollary}[theorem]{Corollary}
\newtheorem{conjecture}[theorem]{Conjecture}

\theoremstyle{definition}

\def\Cee{\mathbb{C}}
\def\Ree{\mathbb{R}}

\def\fH{\mathcal{H}}

\def\fU{\mathcal{U}}

\def\alp{\alpha}
\def\del{\delta}
\def\Del{\Delta}
\def\eps{\varepsilon}

\def\lam{\lambda}

\def\sig{\sigma}
\def\Sig{\Sigma}
\def\vphi{\varphi}

\def\aand{\text{ and }}

\def\iin{\text{ in }}

\def\endpf{{\hfill$\square$\medskip}}
\def\proof{{\noindent{\bf Proof.}\thickspace}}

\def\comp{\raisebox{.2ex}{${\scriptstyle\circ}$}}

\def\spn{\mathrm{span}}
\def\til#1{\tilde{#1}}
\def\wtil#1{\widetilde{#1}}
\def\wbar#1{\overline{#1}}
\def\what#1{\widehat{#1}}
\def\inprod#1#2{\left\langle #1 | #2 \right\rangle}
\def\id{\mathrm{id}}
\def\norm#1{\left\|#1\right\|}

\def\heis{\mathbb{H}}
\def\cent{\mathrm{Z}}
\def\SL{\mathrm{SL}}

\begin{document}

\title[On the algebra $a(G)$]
{
On the subalgebra of a Fourier-Stieltjes algebra generated by pure positive
definite functions
}

\author{Yin-Hei Cheng, Brian E. Forrest and Nico Spronk}

\begin{abstract}
For a locally compact group $G$, the first-named author considered the closed subspace
$a_0(G)$ which is generated by the pure positive definite functions.
In many cases $a_0(G)$ is itself an algebra.  We illustrate using Heisenberg groups and
the $2\times 2$ real special linear group, that this is not the case in general.
We examine the structures of the algebras thereby created and examine properties
related to amenability.  
\end{abstract}

\maketitle

\footnote{{\it Date}: \today.


The second and third named authors were supported by NSERC Discovery Grants.
}

For a locally compact group $G$ let $B(G)$ denote its Fourier-Stieltjes algebra and
$A(G)$ its Fourier algebra, as defined in \cite{eymard}.
The first named author (\cite{cheng}) defined $a_0(G)$ to be the closed linear span
in $B(G)$ of the pure positive definite functions, and then let $a(G)$
denote the closed subalgebra in $B(G)$ generated by $a_0(G)$.  In the case that
$G$ is abelian, and $B(G)=M(\hat{G})$ via Fourier-Stieltjes transform, we have that
$a_0(G)=a(G)\cong\ell^1(\hat{G})$, where the latter is the closed subspace (algebra) generated by
Dirac measures.  Thus we think of the space $a_0(G)$, and the algebra $a(G)$,
as dual analogues of $\ell^1(\hat{G})$.

We use the notation and many results from \cite{arsac}.  We let for a continuous
unitary representation $\pi:G\to\fU(\fH_\pi)$
\[
A_\pi=\wbar{ \spn\{s\mapsto\inprod{\pi(s)\xi}{\eta}:\xi,\eta\in\fH_\pi\} }^{\norm{\cdot}_B}.
\]
We also use the facts that $A_\pi=A_{\pi'}$ if and only if $\pi\simeq\pi'$, i.e.\ the representations
are quasi-equivalent (\cite[(3.1)]{arsac}); 
if $\pi$ and $\sig$ are disjoint, i.e.\ they share no equivalent subrepresentations, 
then $A_{\pi\oplus\sig}=A_\pi\oplus_{\ell^1}A_\sig$ (\cite[(3.13)]{arsac}); 
and $\wbar{\spn A_\pi A_\sig}^{\norm{\cdot}_B}= A_{\pi\otimes\sig}$ (\cite[(3.25)]{arsac}).
We define for any family of representations
$\Sig$, $A_\Sig=\overline{\sum_{\sig\in\Sig}A_\sig}^{\norm{\cdot}_B}$.  
Thus if this family of representaions is pairwise disjoint we have
\[
A_\Sig=\ell^1\text{-}\bigoplus_{\sig\in\Sig}A_\sig.
\]
Thus if $\hat{G}$ denotes the set of irreducible representations, i.e.\
a full set of representatives, one from each unitary equivalence class, we have 
\[
A_{\hat{G}}=a_0(G).
\]
This can be easily seen from the fact that any pure positive definite function
belongs to some $A_\pi$, with $\pi$ in $\what{G}$.
This is thanks to the Gelfand-Naimark-Segal construction characterising pure positive definite
functions, and the fact that each $A_\pi=\spn\{s\mapsto \inprod{\pi(s)\xi}{\xi}:\xi\in\fH_\pi\}$.
If $\hat{G}_F$ denotes the family of finite dimensional irreducible representations
we let
\[
A_F(G)=A_{\hat{G}_F}.
\]
If $\gamma_{ap}:G\to G^{ap}$ denotes the almost periodic compactification, we have
$A_F(G)=A(G^{ap})\comp\gamma_{ap}$ (\cite[(2.97)]{eymard}).

In \cite{cheng} several cases were examined in which $a_0(G)=a(G)$.  This holds
for Moore groups (i.e.\ those groups whose irreducible representations are each finite 
dimensional), in which case $a_0(G)=a(G)=A_F(G)$.  Also, for the $ax+b$-group $G$, we have
$a_0(G)=A(G)\oplus_{\ell^1}A_F(G)=a(G)$.  Our goal is to investigate
$a(G)$ for some cases where $a(G)\not=a_0(G)$, and learn about the structures of these
alegbras.  The two examples considered in the sequel exhibit some features in common, though
quite different structures in terms of operator amenability theory.  See
the survey \cite{spronk1} for more context on amenability properties of  Fourier
and Fourier-Stieltjes algebras.

\section{Heisenberg groups}
We let
\[
\heis_n=\{(p,q,t):p,q\in\Ree^n,t\in\Ree\}
\]
be the Heisenberg group of dimension $2n+1$ with usual matricial group law
\[
(p,q,t)(p',q',t')=(p+p',q+q',t+p\cdot q'+t')
\]
where $p\cdot q$ is the usual dot product of $p$ and $q$.  This is also
called the ``polarized form'' in \cite{folland1}.  Notice that the centre
$\cent$ of $\heis_n$ is given by
\[
\cent=\{(0,0,t):t\in\Ree\}
\]
and there is a natural isomorphism
\[
\heis_n/\cent\cong\Ree^2.
\]
Following \cite{folland} we have that the unitary dual is given by
\begin{equation}\label{eq:heisdual}
\what{\heis}_n=\{\rho_h,\chi_{\xi,\eta}:h\in\Ree^{\not=0},\xi,\eta\in\Ree^n\}
\end{equation}
where the Schr\"{o}dinger representations are given by
\[
\rho_h(p,q,t)f(x)=e^{i(ht+q\cdot x)}f(x+hp),\; f\in L^2(\Ree^2)
\]
(at least up to unitary equivalence)
and the finite dimensional irreducible representations are simply
the characters
\[
\chi_{\xi,\eta}(p,q,t)=e^{i(p\cdot\xi+q\cdot\eta)}.
\]

\begin{proposition}
We have the following tensor product equivalences
\begin{align*}
\chi_{\xi,\eta}\otimes\rho_h&\cong\rho_h \\
\rho_h\otimes\rho_{h'}&\simeq\rho_{h+h'}\text{ if }h+h'\not=0 \\
\rho_h\otimes\rho_{-h}&\cong\lam_{\heis_n/\cent}\comp q
\end{align*}
where $\cong$ is the relation of unitary equivalence and
$\simeq$ that of quasi-equivalence, $\lam_{\heis_n/\cent}$ is the left
regular representation of $\heis_n/\cent$, and $q:\heis_n\to\heis_n/\cent$
is the quotient map.
\end{proposition}

\proof The first two follow, in part, from the Stone-von Neumann theorem:
\begin{equation}\label{eq:stonevonneumann}
\pi\simeq\rho_h\quad\Leftrightarrow\quad
\pi(0,0,t)=e^{iht}I.
\end{equation}
See \cite{folland,folland1}, for example.  Thus unitary equivalence
in the first tensor product follows from the following:  for any group $G$
if we have $\chi\otimes\pi\simeq\pi$ for a character $\chi$ and a representation
$\pi$, then $\chi\otimes\pi\cong\pi$.  Indeed $\{\chi(g)\pi(g):g\in G\}$
can admit only those operators which commute with each $\pi(g)$ as commutators,
and hence by Schur's lemma is irreducible.  Two irreducible representations
are quasi-equivalent only when they admit non-trivial intertwiners and thus are
unitarily equivalent.

Now, let us consider $\rho_h\otimes\rho_{-h}$.  We have
for $f\in L^2(\Ree^n\times\Ree^n)\cong L^2(\Ree^n)\otimes_2L^2(\Ree^n)$ that
\[
\rho_h\otimes\rho_{-h}(p,q,t)f(x,y)
=e^{iq\cdot(x+y)}f(x+hp,y-hp).
\]
We let $W:L^2(\Ree^n\times\Ree^n)\to L^2(\Ree^n\times\Ree^n)$ 
be implemented by
the orthogonal transformation $(x,y)\mapsto\frac{1}{\sqrt{2}}
(x-y,x+y)$ so
\[
Wf(x,y)=f(\textstyle{\frac{1}{\sqrt{2}}}(x-y),\textstyle{\frac{1}{\sqrt{2}}}(x+y))
\text{ and }
W^*f(x,y)=f(\textstyle{\frac{1}{\sqrt{2}}}(x+y),\textstyle{\frac{1}{\sqrt{2}}}(-x+y))
\]
We thus have
\begin{align*}
W\rho_h\otimes\rho_{-h}(p,q,t)W^*f(x,y)
&=[\rho_h\otimes\rho_{-h}(p,q,t)W^*
f](\textstyle{\frac{1}{\sqrt{2}}}(x-y),\textstyle{\frac{1}{\sqrt{2}}}(x+y)) \\
&=e^{i\sqrt{2}q\cdot x}
W^*f((\textstyle{\frac{1}{\sqrt{2}}}(x-y)+hp,\textstyle{\frac{1}{\sqrt{2}}}(x+y)-hp) \\
&=e^{i\sqrt{2}q\cdot x}f(x,y-\sqrt{2}hp) \\
&=(V\otimes I)\lam_{\Ree^n\times\Ree^n}(\sqrt{2}q,\sqrt{2}hp)(V^*\otimes I)f(x,y)
\end{align*}
where $V:L^2(\Ree^n)\to L^2(\Ree^n)$ is the Fourier-Plancherel transform.
Now the map $(p,q,t)\mapsto \sqrt{2}(q,hp):\heis_n\to\Ree^n\times\Ree^n$ is an open 
surjective homomorphism with kernel $\cent$.  The unitary equivalence
$\rho_h\otimes\rho_{-h}\cong\lam_{\heis_n/\cent}\comp q$ follows. \endpf

We are grateful to H.H.~Lee for pointing out the role of the Stone-von Neumann theorem, 
(\ref{eq:stonevonneumann}) above.
K.F.~Taylor kindly informs us that the formula for $\rho_h\otimes\rho_{-h}$
may also be deduced from results of Mackey; details of will be available in his forthcoming
book with E.~Kaniuth.  Our procedure has the benefit of being direct and elementary.  

We let $R=\{\rho_h:h\not=0\}$. From (\ref{eq:heisdual}) we see that
\[
a_0(\heis_n)=A_R\oplus_{\ell^1} A_F(\heis_n).
\]
The conclusions of the proposition above may be reinterpreted as follows:
\begin{align}
\chi_{\xi,\eta}A_{\rho_h}&=A_{\rho_h} \notag \\
\wbar{\spn A_{\rho_h}A_{\rho_{h'}}}^{\norm{\cdot}_B}
&= A_{\rho_{h+h'}}\text{ if }h+h'\not=0 \label{eq:schrodmult} \\
\wbar{\spn A_{\rho_h}A_{\rho_{-h}}}^{\norm{\cdot}_B}
&= A(\heis_n/\cent)\comp q.\label{eq:schrodmult0}
\end{align}
In particular $a_0(\heis_n)$ is not an algebra.

\begin{proposition}\label{prop:algebra}
The closed algebra generated by $a_0(\heis_n)$ is given by
\[
a(\heis_n)=A_R\oplus_{\ell^1}
A(\heis_n/\cent)\comp q\oplus_{\ell^1}A(\heis_n^{ap})\comp\gamma_{ap}.
\]
\end{proposition}

\proof The multiplication relation (\ref{eq:schrodmult0}) shows that $A(\heis_n/\cent)\comp q
\subset a(\heis_n)$.  Each character $\chi_{\xi,\eta}$ clearly multiplies elements of 
$A(\heis_n/\cent)\comp q$, respectively $A_{\rho_h}$, into the same space.
The multiplication relation (\ref{eq:schrodmult}) shows that $A_{R^+}=
\ell^1\text{-}\bigoplus_{h>0}A_{\rho_h}$ and $A_{R^-}=\ell^1\text{-}\bigoplus_{h>0}A_{\rho_h}$
are subalgebras of $a(\heis_n)$.  Finally, we see that $[A(\heis_n/\cent)\comp q]A_{\rho_h}
\subset A_{\rho_h}$.  Indeed, it is immediate from (\ref{eq:stonevonneumann}),
applied to $\lam_{\heis_n/\cent}\comp q\otimes \rho_h$, that this representation
is quasi-equivalent to $\rho_h$.  Hence $\wbar{\spn A_{R^+}A_{R^-}}^{\norm{\cdot}_B}
A_R\oplus_{\ell^1}A(\heis_n/\cent)\comp q$ is an ideal in $a(\heis_n)$.
\endpf

Our goal is to now understand the ideal $A_R\oplus_{\ell^1}A(\heis_n/\cent)\comp q$.
To this end, let us consider a partial compactification of $\heis_n$.  Let
\[
\wtil{\heis}_n=\{(p,q,z):p,q\in\Ree^n,z\in\Ree^{ap}\}.
\]
This group has group law
\[
(p,q,z)(p',q',z')=(p+p',q+q',z\gamma(p\cdot q')z')
\]
where $\gamma:\Ree\to\Ree^{ap}$ is the compactification map, and we write
the group law on $\Ree^{ap}$ multiplicatively.  Let $\til{\gamma}:\heis_n\to
\wtil{\heis}_n$ be the homomorphism given by 
\[
\til{\gamma}(p,q,t)=(p,q,\gamma(t)).
\]

\begin{theorem}\label{theo:partialcompact}
We have
\[
A_R\oplus_{\ell^1}A(\heis_n/\cent)\comp q
=A(\wtil{\heis}_n)\comp\til{\gamma}.
\]
In particular, this algebra has Gelfand spectrum
isomorphic to $\wtil{\heis}_n$.
\end{theorem}

\proof 
We begin by noting that the Haar measure on $\wtil{\heis}_n$ is the product measure
$m_n\times m_n\times\mu$ where $m_n$ is the Lebesgue measure and
$\mu$ is the Haar measure on $\Ree^{ap}$; indeed
\[
\iint_{\Ree^{2n}}\int_{\Ree^{ap}}\vphi(p+p',q+q',z\gamma(p\cdot q')z')\,d(p,q)\,dz
=\iint_{\Ree^{2n}}\int_{\Ree^{ap}}\vphi(p,q,z)\,d(p,q)\,dz
\]
for compactly supported continuous $\vphi:\wtil{\heis}_n\to\Cee$.
Thus we obtain a decomposition
\[
L^2(\wtil{\heis}_n)\cong L^2(\Ree^{2n})\otimes_2 \left(\ell^2\text{-}
\bigoplus_{h\in\Ree}\Cee\chi_h\right)
\]
where $\chi_h$ is the character on $\Ree^{ap}$ given by $h\iin\Ree_d\cong\what{\Ree^{ap}}$.
We let 
\[
\fH_h=L^2(\Ree^{2n})\otimes_2\Cee\chi_h\cong L^2(\Ree^{2n})
\]
so for $f\iin\fH_h$, $f(x,y,zz')=\chi_h(z)f(x,y,z')$.
We compute for such $f$ the left regular representation
\begin{align*}
\lam_{\wtil{\heis}_n}(p,q,z)f(x,y,z')&=f((-p,-q,\gamma(p\cdot q)z^{-1})(x,y,z')) \\
&=f(x-p,y-q,\gamma(p\cdot(q-y))z^{-1}z') \\
&=e^{ihp\cdot(q-y)}\chi_{-h}(z)f(x-p,y-q,z')
\end{align*}
We immediately observe that
\[
\lam_{\wtil{\heis}_n}(\cdot)|_{\fH_h}\comp\til{\gamma}(0,0,t)=\chi_{-h}\comp\gamma(t)I=e^{-iht}I.
\]
Thus it follows (\ref{eq:stonevonneumann}) that $\lam_{\wtil{\heis}_n}(\cdot)|_{\fH_h}\comp\til{\gamma}
\simeq\rho_{-h}$ for $h\not=0$.  We also immediately see that
$\lam_{\wtil{\heis}_n}|_{\fH_0}\comp\til{\gamma}\cong\lam_{\heis_n/\cent}\comp q$.
(Observe that $\fH_0$ is the only component of $L^2(\wtil{\heis}_n)$ which is not 
annihilated by averaging over the centre $\wtil{\cent}$ of $\wtil{\heis}_n$.)

The identification of the algebras is immediate.  The identification
of the spectrum follows \cite[(3.34)]{eymard}.  \endpf

Any irreducible representation of $\wtil{\heis}_n$ must also be an irreducible representation
of $\heis_n$.  Hence it is immedate that $a(\heis_n)=a(\wtil{\heis}_n)\comp \til{\gamma}\cong
a(\wtil{\heis}_n)$.  This is despite that $\heis_n\not\cong\wtil{\heis}_n$.  Of course, a similar
phenomenon may be observed with any non-compact Moore group $G$:
$a(G)=A(G^{ap})\circ\gamma_{ap}\cong a(G^{ap})$.

The spine $A^*(G)$ for a locally compact group $G$ is defined in \cite{ilies}.
It may be given as
\[
A^*(G)=\wbar{\sum_{(\eta,H)}A(H)\comp\eta}
\]
where the sum runs over all pairs where $\eta:G\to H$ is a continuous homomorphism
into a locally compact group with dense range.  Details as to why this sum can be determined
as a sum over an index {\em set} are given in \cite{ilies}.

\begin{corollary}\label{cor:amenable}
{\bf (i)} We have $a(\heis_n)\subset A^*(\heis_n)$.

{\bf (ii)} The algebra $a(\heis_n)$ is operator amenable, but not amenable.
\end{corollary}

\proof We have $a(\heis_n)=A(\wtil{\heis}_n)\comp \til{\gamma}\oplus_{\ell^1}A(\heis_n^{ap})
\comp\gamma_{ap}$, which clearly gives (i).  That $a(\heis_n)$ is operator amenable
is an immediate consequence of \cite[Prop.\ 3.1]{rundes}.  We observe that
$a(\heis_n)$ admits $a(\wtil{\heis}_n)\comp \til{\gamma}\cong
a(\wtil{\heis}_n)$ as a complemented ideal.  This ideal is not amenable thanks to
\cite[Thm.\ 2.3]{forrestr} of \cite[Cor.\ 3.3]{runde}.  Thus by \cite[2.3.7]{runde},
$a(\heis_n)$ is not amenable.
\endpf

Motivated by all of the examples we have thus far, we are emboldened
to suggest the following.  We let $vn(G)=\ell^\infty\text{-}\bigoplus_{\pi\in\hat{G}}\mathcal{B}(\fH_\pi)$,
which is a von Neumann algebra and the dual of $a(G)$. We also refer to
\cite{berglund,ilies,sstokke} for information on topological Clifford semigroups.

\begin{conjecture}\label{conj:spine}
{\bf (i)} For any locally compact group $G$, there is an injective continuous homomorphism
into a locally compact group with dense range, $\gamma:G\to H$, such that
$A(H)\comp\gamma$ is contained in $a(G)$ and is an ideal in $a(G)$.

{\bf (ii)} The Gelfand spectrum $\Phi_{a(G)}\subset vn(G)$ is a Clifford semigroup
with a dense open subgroup isomorphic to $H$.
\end{conjecture}

Indeed, for Moore groups we use the almost periodic compactification
$\gamma_{ap}:G\to G^{ap}$, and $\Phi_{a(G)}\cong G^{ap}$.  
For $G$ being either of the $ax+b$-group, $\SL_2(\Ree)$
(see below), or $\wtil{\heis}_n$, we use $\id:G\to G$;
and $\Phi_{a(G)}=G\sqcup G^{ap}$.  For $G=\heis_n$ we use
$\til{\gamma}:\heis_n\to\wtil{\heis}_n$, and $\Phi_{a(\heis_n)}=\heis_n\sqcup\wtil{\heis}_n\sqcup
\heis_n^{ap}$.   The truth of (i), above, would verify a conjecture in \cite{cheng}, that
the invertible part of $\Phi_{a(G)}$ consists of unitaries.

Of course for a discrete non-Moore group, i.e.\ not Type I (see \cite[12.6.37]{palmer}),
we will not be able to calculate $a_0(G)$, nor $a(G)$, in the elementary manner presented here.

\section{$\SL_2(\Ree)$}

We show how results of Repka \cite{repka}, on the tensor products of irreducible
representations on $\SL_2(\Ree)$, give a structure theory for $a(\SL_2(\Ree))$.

We begin by listing all of the families of irreducible unitary representations.  Our notation
is similar to that of \cite[p.\ 247]{folland}, with the exception of our parametrisation of
the complementary series.  We shall, not need, and thus will not indicate, any of the 
actual formulas for the representations themselves.
\begin{flalign*}
&\text{principal continuous series }&\Pi^+&=\{\pi_t^+:t\in[0,\infty)\},\, \Pi^-=\{\pi_t^-:t\in(0,\infty)\} \\
&\text{complementary series }&K&=\{\kappa_s:s\in(-1,0)\} \\
&\text{discrete series }&\Del^{\pm}&=\{\del^\pm_n:n=2,3,4,\dots\} \\
&\text{mock discrete series }&M&=\{\del_1^+,\del^-_1\}. 
\end{flalign*}
There is also the trivial representation $1$. 
We will consider two direct integral representations and a direct sum:
\[
\pi^\pm=\int_{(0,\infty)}^\oplus\pi_t^{\pm}\,dt,\quad
\del=\bigoplus_{n=2}^\infty(\del_n^+\oplus\del_n^-).
\]
The trace formula of Harish-Chandra (\cite{harishchandra}) tells us that there is a quasi-equivalence
$\lam\simeq\pi^+\oplus\pi^-\oplus\del$.
In other words 
\[
A(G)=A_{\pi^+}\oplus_{\ell^1}A_{\pi^-}\oplus_{\ell^1}A_\del.
\]

We record a crude summary of \cite{repka}, Theorems
4.6, 5.9, 6.4, 7.1, 7.3 and 8.1.

\begin{lemma}\label{lem:repka}
Let $\sig,\tau$ be any two non-trivial 
irreducible unitary representations of $\SL_2(\Ree)$. 
Then we have quasi-containments
\[
\sig\otimes\tau\prec\begin{cases}
\pi^+\oplus\pi^-\oplus\del\oplus\kappa_{s+t+1}&\text{if }\{\sig,\tau\}=\{\kappa_s,\kappa_t\}
\aand s+t<-1 \\  \pi^+\oplus\pi^-\oplus\del &\text{otherwise.}\end{cases}
\]
\end{lemma}

\begin{theorem}\label{theo:sltwor}
{\bf (i)} We have 
\begin{align*}
a(\SL_2(\Ree))&=a_0(\SL_2(\Ree))\oplus_{\ell^1}A_{\pi^+}\oplus_{\ell^1}A_{\pi^-} \\
&=A_{\Pi^+}\oplus_{\ell^1}A_{\Pi^-}\oplus_{\ell^1}A_K\oplus A_M
\oplus_{\ell^1}\Cee1\oplus_{\ell^1}A(G).
\end{align*}
{\bf (ii)} We have
\[
\wbar{\spn a(\SL_2(\Ree))^2}=A_K\oplus_{\ell^1}\Cee1\oplus_{\ell^1}A(G).
\]
\end{theorem}

\proof Clearly $a_0(\SL_2(\Ree))=A_{\Pi^+}\oplus_{\ell^1}A_{\Pi^-}\oplus_{\ell^1}A_K\oplus A_M
\oplus_{\ell^1}\Cee1\oplus_{\ell^1} A_{\Del^+}\oplus_{\ell^1}A_{\Del^-}$.    For $\Sig,T$ being any of
$\Pi^\pm,M,\Del^\pm$, Lemma \ref{lem:repka} tells us that
\[
A_\Sig A_K,\,A_\Sig A_T\subset A_{\pi^+}\oplus_{\ell^1}A_{\pi^-}\oplus_{\ell^1}A_\del=A(G).
\]
whereas
\[
A_K^2\subset A_{\pi^+}\oplus_{\ell^1}A_{\pi^-}\oplus_{\ell^1}A_\del\oplus_{\ell^1}A_K=
A(G)\oplus_{\ell^1}A_K.
\]
Hence both (i) and (ii) are immediate. \endpf

\begin{corollary}
The Gelfand spectrum of $a(\SL_2(\Ree))$ is the one-point compactification,
$\SL_2(\Ree)_\infty$.
\end{corollary}

\proof  We first observe that  $\SL_2(\Ree)_\infty$ is the spectrum
of $A(\SL_2(\Ree))\oplus\Cee1$.  It can be easily derived from Lemma \ref{lem:repka}
that if $u=\sum_{k=1}^n u_k$ where $u_k\in A_{\pi_k}$ for $\pi_k\in\Pi^+\cup\Pi^-\cup M\cup K$
then $u^n\in A(\SL_2(\Ree))$ for some $n$.  Indeed, if, up to reordering of indices $k$,
we have that $\pi_k=\kappa_{s_k}$, $s_1<s_2\dots<s_l$ for some $l\leq m$, $\pi_k\not\in K$ for 
$k>l$, then $n\leq\log_2(\frac{1}{1+s_1})-1$.  We have $n=2$, otherwise.
Hence it follows that for $u\in a(\SL_2(\Ree))$ and $\eps>0$ that there is $v$ in $a(\SL_2(\Ree))$
for which $\|u-v\|<\eps$ and $v^n\in A(G)\oplus\Cee 1$.  \endpf

The corollary above can also be deduced from the main result of \cite{cowling}, and has
a similar method of proof.

We observe that $a(\SL_2(\Ree))$ admits much weaker amenability properties
than does $a(\heis_n)$.

\begin{corollary}
The algebra $a(\SL_2(\Ree))$ admits no non-zero point derivations, but 
is not operator weakly amenable.
\end{corollary}

\proof  Since $A(\SL_2(\Ree))$ is operator weakly amenable (\cite{spronk}),
it admits no non-zero point derivations.  For $u$ in $a(\SL_2(\Ree))$ and $\eps>0$, the
proof of the corollary
above provides $\alp\in\Cee$ and $v,w$ in $a(\SL_2(\Ree))$ for which $u=v+w+\alp 1$ where 
$\|w\|\leq\eps$ and $u^n\in A(\SL_2(\Ree))$ for some $n$.  It follows that
any point derivation is zero.  Theorem \ref{theo:sltwor} (ii) shows that
$\spn a(\SL_2(\Ree))^2$ is not dense in $a(\SL_2(\Ree))$, hence $a(\SL_2(\Ree))$
is not operator weakly amenable by \cite[Lem.\ 3.1]{forrestw}. \endpf

In \cite{ilies} it is computed that $A^*(\SL_2(\Ree))=A(\SL_2(\Ree))\oplus_{\ell^1}\Cee 1$.
We let, for any locally compact group $G$, the Rajchman algebra $B_0(G)$ be the subalgebra of
$B(G)$ consisting of those elements vanishing at $\infty$.  It is observed in \cite{cowling}
that $B(\SL_2(\Ree))=B_0(\SL_2(\Ree))\oplus_{\ell^1}\Cee 1$.  For any continuous
singular Borel measure $\mu$ on $(0,\infty)$ we have that $\pi_\mu^\pm=\int^\oplus_{(0,\infty)}
\pi_t^{\pm}\,d\mu(t)$ satisfies $A_{\pi_\mu^\pm}\cap a(\SL_2(\Ree))=\{0\}$; indeed
see \cite[(3.12) \& (3.55)]{arsac}.  Hence
\[
A^*(\SL_2(\Ree))\subsetneq a(\SL_2(\Ree))\subsetneq B_0(\SL_2(\Ree))\oplus_{\ell^1}\Cee 1
=B(\SL_2(\Ree)).
\]


\end{document}